\let\OToneAccents=\relax
\magnification=\magstep1
\frenchspacing
\baselineskip=16truept
\font\douze=cmr10 at 12pt

\font\grande=cmb10 at 16truept
\def\titre#1{{\OToneAccents\noindent\grande #1}}

\def \chapter#1{\vfill\eject\ifodd\pageno\else\ \vfill\eject\fi\centerline{\grande #1}\bigskip}

\font\tendo=wncyr10 at 12pt
\font\sevendo=wncyr7
\newfam\dofam 
\textfont\dofam = \tendo 
\scriptfont\dofam= \sevendo

\font\tendb=msbm10 
\font\sevendb=msbm7
\newfam\dbfam 
\textfont\dbfam = \tendb 
\scriptfont\dbfam= \sevendb
\def\db {\fam\dbfam\tendb}
\font\tenrsfs=rsfs10 
\font\sevenrsfs=rsfs7
\newfam\scrfam 
\textfont\scrfam = \tenrsfs
\scriptfont\scrfam= \sevenrsfs

\def\F{{\bi F}}

\def\Q{{\db Q}}

\def\Z{{\db Z}}

\def\F{{\db F}}

\def\Inf{\mathop{\rm Inf}\nolimits}

\font\tendb=msbm10 
\font\sevendb=msbm7
\newfam\dbfam 
\textfont\dbfam = \tendb 
\scriptfont\dbfam= \sevendb
\def\db {\fam\dbfam\tendb}
\font\tenrsfs=rsfs10 
\font\sevenrsfs=rsfs7
\newfam\scrfam 
\textfont\scrfam = \tenrsfs
\scriptfont\scrfam= \sevenrsfs

\def\F{{\bi F}}

\def\Q{{\db Q}}
\def\Z{{\db Z}}
\def\N{{\db N}}

\def\F{{\db F}}

\def\mod{\hbox{\rm \ mod.\ }}

\font\teufm=eufm10
\font\seufm=eufm10 at 7pt
\font\sseufm=eufm10 at 6pt
\newfam\fameufm
\textfont\fameufm=\teufm
\scriptfont\fameufm=\seufm
\scriptscriptfont\fameufm=\sseufm
\def\goth{\fam\fameufm\teufm}

\def\wp{{\goth p}}

\def\qq{{\goth q}}

\def\mod{\hbox{\rm \ mod.\ }}

\def\og{\leavevmode\raise.3ex\hbox{$\scriptscriptstyle\langle\!\langle\,$}}
\def\fg{\leavevmode\raise.3ex\hbox{$\scriptscriptstyle\,\rangle\!\rangle$}}

\centerline {\titre{Remarques sur le   premier cas du th\'eor\`eme de Fermat}}
\medskip
\centerline {\titre{ sur les corps de nombres}}
\bigskip
\smallskip
\centerline {Alain Kraus}
\bigskip 
\smallskip
{\bf{Abstract.}} The first case of Fermat's Last Theorem    for a prime  exponent $p$ can sometimes be proved  using the existence of local  obstructions.
In 1823,  Sophie Germain has obtained an important  result in this direction by establishing that, if $2p+1$ is a prime number,   the first case of Fermat's Last Theorem is true  for  $p$. 
In this paper, we investigate  such obstructions   over number fields. 
We  obtain   analogous results on Sophie Germain type criteria,   for  imaginary quadratic fields. Furthermore,  extending a well known statement over $\Q$, we give 
an easily testable condition which allows occasionally  to prove  the first  case of Fermat's Last Theorem over  number fields  for a prime number  $p\equiv 2 \mod 3$.
\bigskip

\medskip

\bigskip
\medskip

\centerline{INTRODUCTION}
\bigskip
Soient $K$ un corps de nombres et $p$ un nombre premier impair. On dit que le premier cas du  th\'eor\`eme de Fermat est vrai sur $K$ pour 
l'exposant $p$, s'il n'existe pas  d'\'el\'ements $x,y,z$ dans l'anneau d'entiers $O_K$ de $K$ tels que  
$$x^p+y^p+z^p=0\quad \hbox{et}\quad  (xyz)O_K+pO_K=O_K.$$
En 1994, A. Wiles a d\'emontr\'e qu'il en est ainsi pour le corps $\Q$, ind\'ependamment du fait que $xyz$ soit premier avec $p$ ([Wi]). En 2004, ce r\'esultat a \'et\'e \'etendu au corps $\Q\bigl(\sqrt{2}\bigr)$  par F. Jarvis et P. Meekin ([Ja-Me]). R\'ecemment, N. Freitas et S. Siksek ont accompli 
des progr\`es importants concernant  l'\'equation de Fermat sur  les corps  totalement r\'eels 
 ([Fr-Si]).
\smallskip
On s'int\'eresse dans cet article au premier cas du th\'eor\`eme de Fermat. Il est bien connu que la pr\'esence d'obstructions locales  permet  de le d\'emontrer sur $\Q$ pour certains exposants $p$. 
Historiquement, le premier r\'esultat spectaculaire  \`a ce sujet  a \'et\'e obtenu par  S. Germain en 1823, qui a d\'emontr\'e    que  si $2p+1$ est un nombre premier,  le premier cas du th\'eor\`eme de Fermat est vrai sur $\Q$ pour l'exposant $p$.  Cet \'enonc\'e a \'et\'e g\'en\'eralis\'e par de nombreux math\'ematiciens, notamment 
par Wendt  en 1884 qui a obtenu une g\'en\'eralisation en termes du r\'esultant des polyn\^omes de la forme  $X^n-1$ et $(X+1)^n-1$ (voir [Co-2] p. 430 et [Ri] p. 137). 
On se propose  ici de d\'emontrer un analogue du crit\`ere  de Wendt  dans le cas o\`u $K$ est un corps quadratique imaginaire. 
\vskip0pt\noindent
Par ailleurs,   si pour tout entier $a$ compris entre $1$ et ${p-3\over 2}$,  on a 
$$1+a^p\not\equiv (1+a)^p \mod p^2,$$
le premier cas du th\'eor\`eme de Fermat est vrai sur $\Q$ pour l'exposant $p$ ([Co-2] p. 430). 
Cette condition ne peut \^etre satisfaite que si $3$ ne divise pas $p-1$. F.H. Hao et C.J. Parry ont \'etendu ce crit\`ere aux corps quadratiques  dont l'anneau d'entiers poss\`ede un id\'eal premier au-dessus de $p$ de degr\'e r\'esiduel 1 ([Ha-Pa]).
On obtient  ici  une g\'en\'eralisation de cet \'enonc\'e \`a d'autres familles de corps de nombres. Pour tout  nombre premier $p\geq 5$  v\'erifiant cette  condition modulo $p^2$, cela permet par exemple  d'\'etablir le premier cas du th\'eor\`eme de Fermat pour $p$, sur les corps cubiques purs, et sur des corps purs  de degr\'e sur $\Q$ arbitrairement grand.  
\medskip
J'ai b\'en\'efici\'e de nombreuses remarques de D. Bernardi  pendant la r\'edaction de  cet article. Je l'en remercie vivement ici.
\bigskip

{\douze I. \'Enonc\'e des r\'esultats}
\medskip
Pour tout corps de nombres $K$, notons  $O_K$ son anneau d'entiers et $h_K$ son nombre de classes. La lettre $p$ d\'esigne un nombre premier impair. 
\bigskip

{\bf{1. Le crit\`ere de Wendt sur les corps quadratiques imaginaires}}
\medskip
Soit $K$ un corps quadratique imaginaire.
Pour tout entier $n\geq 1$, notons $W_n$ le r\'esultant des polyn\^omes  $X^n-1$ et $(X+1)^n-1$.
\bigskip

\proclaim Th\'eor\`eme 1.  Supposons   les conditions suivantes  satisfaites :
\smallskip
\vskip0pt\noindent
1) on a $h_K\not \equiv 0 \mod p$.
\smallskip
\vskip0pt\noindent
2) Il existe  $n\geq 1$ tel que  $q=np+1$ soit un nombre premier d\'ecompos\'e dans $K$  et que 
$$(n^n-1)W_n\not\equiv 0 \mod q.$$
\vskip0pt\noindent
Alors, le premier cas du th\'eor\`eme de Fermat est vrai sur $K$ pour l'exposant $p$.
\bigskip

{\bf{Remarque 1.}} Le nombre premier $p$ \'etant donn\'e, 
il  n'existe
 qu'un nombre fini d'entiers $n$ pour lesquels  $np+1$ soit un nombre premier ne divisant pas $W_n$.  Plus pr\'ecis\'ement, 
 Dickson  a d\'emontr\'e en 1909 que 
tout nombre premier de la forme $np+1$, plus grand que
 $(p-1)^2(p-2)^2+6p-2,$
 divise $W_n$ ([Ri] p. 301). De plus,   $W_n=0$ si  $6$ divise $n$. La  question de savoir  si pour tout  $p$, il existe  $n$ tel que $np+1$ soit un nombre premier   ne divisant  pas $W_n$ est toujours ouverte. Elle  a \'et\'e pos\'ee  par Flye Sainte-Marie en 1880 et par Landau en 1913 ({\it{loc. cit.}}). 
Pour autant, on constate exp\'erimentalement que la seconde condition du th\'eor\`eme est tr\`es souvent  r\'ealis\'ee en pratique, except\'e comme il se doit pour le corps $\Q\big(\sqrt{-3}\bigr)$  o\`u elle ne l'est jamais si $p\geq 5$, \`a cause de la pr\'esence  des racines cubiques de l'unit\'e. Si $K$ n'est pas $\Q\big(\sqrt{-3}\bigr)$, il est plausible  qu'elle le soit  toujours   d\`es que  $p$ est plus grand qu'une constante ne d\'ependant que de $K$. 
\bigskip
On d\'eduit de l'\'egalit\'e $W_2=-3$ un analogue du r\'esultat de S. Germain sur $K$.
\bigskip

\proclaim Corollaire 1. Supposons $h_K$ non divisible par $p$ et que $2p+1$ soit un nombre premier d\'ecompos\'e dans $K$. Alors, le premier cas du th\'eor\`eme de Fermat est vrai sur $K$ pour  l'exposant  $p$. 
\bigskip

\proclaim Corollaire 2. 1) Si l'on a $p\leq 10^6$, le premier cas du th\'eor\`eme de Fermat est vrai sur le corps $\Q(i)$ pour l'exposant $p$.
\smallskip
\vskip0pt\noindent
2) Si  $4p+1$  ou $8p+1$ ou $16p+1$ est  premier,  le premier cas du th\'eor\`eme de Fermat est vrai sur $\Q(i)$ pour l'exposant $p$.
 \bigskip
 
D\'emonstration :  On a $h_{\Q(i)}=1$ et l'on  peut v\'erifier  la seconde condition du th\'eor\`eme  pour $p\leq 10^6$ \`a l'aide du  logiciel de calculs Pari, environ en cinq minutes ([Pari]).  Tout nombre premier congru \`a $1$ modulo $4$ est d\'ecompos\'e dans $\Q(i)$. Les factorisations de $(n^n-1)W_n$ pour $n=4,8,16$, 
permettent alors d'\'etablir la seconde assertion (cf. {\it{loc. cit.}}). 
\medskip
 \`A titre indicatif, sur le corps $\Q(i)$, la liste des couples $(p,n)$ pour $p<100$, 
 avec les plus petits entiers $n$ pour lesquels le crit\`ere fonctionne est la suivante :
$$(3,4), (5,8), (7,4), (11,8), 13,4), (17,8), (19,40), (23,20), (29,8), (31,76), (37,4), (41,20),$$
$$(43,4), (47,20), (53,20), (59,20), (61,16), (67,4), (71,8), (73,4), (79,4), $$
$$(83,32), (89,44), (97,4).$$
\smallskip

 {\bf{2. Obstructions locales modulo $p^2$}}
 \medskip

\proclaim Th\'eor\`eme 2.  Soient  $K$ un corps de nombres et $p$ un nombre premier impair. Supposons   les  conditions suivantes satisfaites~:
\smallskip
\vskip0pt\noindent
1)  il  existe  un id\'eal premier  de $O_K$ au-dessus de $p$, de degr\'e r\'esiduel  sur $p$ \'egal \`a $1$, et d'indice de ramification sur $p$ inf\'erieur ou \'egal \`a $p-1$. 
\smallskip
\vskip0pt\noindent
2)  On a 
$$1+a^p\not\equiv (1+a)^p \mod p^2 \quad \hbox{pour tout} \ \  a=1,2,\cdots,{p-3\over 2}.\leqno(1)$$
Alors, le premier cas du th\'eor\`eme de Fermat est vrai  sur $K$ pour l'exposant $p$.
\bigskip

{\bf{Remarque 2.}}  La condition (1) de l'\'enonc\'e ne peut \^etre r\'ealis\'ee que si $p=3$ ou si $p\equiv 2 \mod 3$. En effet, si $p\neq 3$, le polyn\^ome 
$(X+1)^p-X^p-1$ est divisible par $p(X^2+X+1)$ et
les racines cubiques de l'unit\'e appartiennent \`a $\F_p$ si $p\equiv 1 \mod 3$.
L'ensemble des nombres premiers $p<$ 150 qui la v\'erifie  est 
$$\Big\lbrace 3, 5, 11, 17, 23, 29, 41, 47, 53, 71, 89, 101, 107, 113, 131, 137, 149\Big\rbrace.$$
Exp\'erimentalement, on constate qu'environ  84 pour cent  des nombres premiers congrus \`a $2$ modulo $3$ satisfont la condition (1). Par exemple, il y a 39265  nombres premiers impairs, congrus \`a $2$ modulo $3$,   plus petits que $10^6$, et  
$33316$  d'entre eux passent  positivement le test ;  il y en a  30870 qui sont  irr\'eguliers et  13192 d'entre eux  satisfont la condition (1).
\bigskip

\proclaim Corollaire 3. Soit $p$ un nombre premier $\geq 5$ v\'erifiant la condition $(1)$. Soit $d$ un entier rationnel distinct de $\pm 1$. Supposons que l'on soit dans l'un des cas  suivants :
\smallskip
\vskip0pt\noindent
1) $K=\Q\bigl(\root 3\of {d}\bigr)$
o\`u $d$ est  sans facteurs cubiques, 
\medskip
\vskip0pt\noindent
2) $K=\Q\bigl(\root n\of {d}\bigr)$ o\`u $p$ ne divise pas $dn$,  $d$ est sans facteurs carr\'es et $n\equiv 1 \mod p-1$,
\smallskip
\vskip0pt\noindent
3) $K$ est une extension de $\Q$, totalement ramifi\'ee en $p$, dont le  degr\'e sur $\Q$ est inf\'erieur ou \'egal \`a  $p-1$.
\smallskip
\vskip0pt\noindent
Alors, le premier cas du th\'eor\`eme de Fermat est vrai sur $K$ pour l'exposant $p$.
 \bigskip

{\bf{Remarque 3.}} Dans l'\'enonc\'e de la  seconde assertion, l'hypoth\`ese selon laquelle $d$ est distinct de $\pm 1$ et sans facteurs carr\'es, sert  \`a garantir  que  le polyn\^ome $X^n-d$ est irr\'eductible sur $\Q$. Par ailleurs, 
si l'on sp\'ecifie $d$ et $p$, on peut  obtenir un \'enonc\'e plus pr\'ecis. \`A titre indicatif, si  $n$ est un entier impair non multiple de $5$, le  premier cas du th\'eor\`eme de Fermat est vrai  sur le corps $\Q\bigl(\root n\of {3}\bigr)$,  pour l'exposant $p=5$.

\bigskip
{\douze II. D\'emonstration du th\'eor\`eme 1}
\medskip
Soit   $(x,y,z)$ un triplet d'\'el\'ements de $O_K$ tel que 
$$x^p+y^p+z^p=0 \quad \hbox{et}\quad (xyz)O_K+pO_K=O_K.$$  
Posons 
$$D=xO_K+yO_K.$$
On a les \'egalit\'es
$$D=xO_K+zO_K=yO_K+zO_K.\leqno(2)$$
\medskip

{\bf{1. Lemmes pr\'eliminaires}}
\medskip
Les \'egalit\'es (2) et le  lemme qui suit n'utilisent pas le fait  que $K$ est un corps  quadratique imaginaire. 
\medskip

\proclaim Lemme 1. L'id\'eal 
$${(x+y)O_K\over D}$$
est la puissance $p$-i\`eme d'un id\'eal de $O_K.$
\bigskip

D\'emonstration :  On a l'\'egalit\'e
$$(x+y) s=-z^p\quad \hbox{o\`u}\quad s=\sum_{k=0}^{p-1} x^{p-1-k}(-y)^k.$$
Parce que $x$ et $y$ sont dans $D$, l'id\'eal $D^{p-1}$ divise $sO_K$. 
On obtient l'\'egalit\'e d'id\'eaux de~$O_K$
$$\biggl({(x+y)O_K\over D}\biggr) \biggl({sO_K\over D^{p-1}}\biggr)=\biggl({zO_K\over D}\biggr)^p.$$
Il suffit ainsi d'\'etablir  que l'on a
$${(x+y)O_K\over D}+ {sO_K\over D^{p-1}}=O_K.$$
Soit $\qq$ un id\'eal premier non nul de $O_K$ divisant ${(x+y)O_K\over D}$. Il s'agit de montrer que $\qq$ ne divise pas 
$ {sO_K\over D^{p-1}}$. Pour cela, on v\'erifie  par r\'ecurrence 
que pour tout $k\geq 1$, on a
$$(-y)^k\equiv x^k \mod \qq D^k.$$
Pour tout $k$ compris entre $0$ et $p-1$, on a donc 
$$x^{p-1-k}(-y)^k\equiv x^{p-1} \mod \qq D^{p-1},$$
d'o\`u la congruence
$$s\equiv px^{p-1}\mod \qq D^{p-1}.$$
Supposons  que  $\qq$ divise ${sO_K\over D^{p-1}}$. L'id\'eal  $\qq D^{p-1}$ divise alors $(px^{p-1})O_K$. Par ailleurs, l'\'egalit\'e $(xyz)O_K+pO_K=O_K$
entra\^\i ne que 
$\qq D^{p-1}$  est premier avec $pO_K$, 
donc $\qq D^{p-1}$ divise $x^{p-1}O_K$. On en d\'eduit que $\qq$ divise
${xO_K\over D}$, puis que $x$ est dans $\qq D$. L'\'el\'ement  $x+y$ \'etant aussi dans $\qq D$, il en est de m\^eme de $y$.
Cela contredit le fait que $D$ est le plus grand commun diviseur  de $xO_K$ et $yO_K$, d'o\`u le lemme.
\bigskip

Le fait que $K$  soit un corps quadratique imaginaire intervient d\'esormais de fa\c con essentielle. On supposera de plus, ce qui n'est pas restrictif,
$$K\neq \Q\bigl(\sqrt{-3}\bigr).$$
En effet, on a $h_{\Q\bigl(\sqrt{-3}\bigr)}=1$, la seconde  condition du th\'eor\`eme 1 est  satisfaite pour  $p=3$ (avec $n=2$) et ne l'est pas si $p\geq 5$.
Par ailleurs,  il  est connu que le premier cas du th\'eor\`eme de Fermat est vrai sur $\Q\bigl(\sqrt{-3}\bigr)$ pour l'exposant $p=3$. 
Le th\'eor\`eme 1 est donc vrai pour le corps $\Q\bigl(\sqrt{-3}\bigr)$.
\bigskip

Par hypoth\`ese, $p$ ne divise pas $h_K$. Il existe  donc $t\in \N$ tel que $p$ divise $th_K+1$.
L'id\'eal $D^{h_K}$ est principal. En particulier, il existe  $d\in O_K$ tel que l'on ait 
 $$D^{h_Kt}=dO_K.\leqno(3)$$
 \medskip
 
\proclaim Lemme 2. Il existe des \'el\'ements non nuls $a,b,c$ dans $O_K$ tels que l'on ait
$$d(x+y)=a^p,\quad  d(x+z)=b^p,\quad d(y+z)=c^p.$$

D\'emonstration :  D'apr\`es le lemme 1, il existe un id\'eal $I$ de $O_K$ tel que  l'on ait l'\'egalit\'e  $(x+y)O_K=DI^p$. On a donc
 $$\Bigl(D^{{h_Kt+1\over p}}I\Bigr)^p=d(x+y)O_K.$$
Parce que $p$ ne divise pas $h_K$, l'id\'eal $D^{{h_Kt+1\over p}}I$ est principal. 
Le corps $K$ \'etant quadratique imaginaire, distinct de $\Q\bigl(\sqrt{-3}\bigr)$, 
les unit\'es de $O_K$ sont des puissances $p$-i\`emes dans $O_K$ (y compris si $p=3$), d'o\`u  l'existence d'un \'el\'ement $a\in O_K$ tel que 
$d(x+y)=a^p.$ 
Les \'egalit\'es (2) et  (3) entra\^\i nent alors le r\'esultat.
\bigskip

Soit $n\geq 1$ un entier tel que $q=np+1$ soit un nombre premier v\'erifiant la seconde condition de l'\'enonc\'e du th\'eor\`eme.
\bigskip

\proclaim Lemme 3. Chaque id\'eal premier de $O_K$ au-dessus de $q$ divise $(xyz)O_K$.

D\'emonstration: Supposons qu'il existe un id\'eal premier $\qq$ de $O_K$ au-dessus de $q$ ne divisant pas $(xyz)O_K$.  Il  existe alors $u\in O_K$ tel que l'on ait
$$u\equiv \Bigl({x\over z}\Bigr)^{{q-1\over n}} \mod \qq.$$
Le corps  $O_K/\qq$ est  de cardinal $q$, d'o\`u la congruence
$$u^{n}= 1 \mod \qq.$$
L'\'egalit\'e $x^p+y^p+z^p=0$ implique 
$$u+1\equiv -\Bigl({y\over z}\Bigr)^{{q-1\over n}} \mod \qq.$$
On obtient ($n$ est pair)
$$(u+1)^n\equiv 1 \mod \qq,$$
ce qui entra\^\i ne  que $q$ divise $W_n$, d'o\`u une contradiction et le r\'esultat.
\bigskip

{\bf{2. Fin de la d\'emonstration  du th\'eor\`eme 1}} 
\medskip
 Quitte \`a  diviser l'\'egalit\'e $x^p+y^p+z^p=0$ par une puissance convenable de $q$,  on peut  supposer que $(x,y,z)$ n'est pas nul modulo $qO_K$. 
Il existe donc un id\'eal premier $\qq$ de $O_K$ au-dessus de $q$ tel que $(x,y,z)$ soit non nul modulo $\qq$. Notons $v_{\qq}$ la valuation sur $K$ qui lui est associ\'ee. 
 D'apr\`es le lemme 3, on peut supposer que l'on a 
 $$v_{\qq}(z)\geq 1,$$
auquel cas on a 
$$v_{\qq}(x)=v_{\qq}(y)=0.$$
En particulier, $\qq$ ne divise pas $D$ et d'apr\`es l'\'egalit\'e (3), on a 
$$v_{\qq}(d)=0.$$
Par suite, on a $v_{\qq}\bigl(d(y+z)\bigr)=0$ et d'apr\`es le lemme 2 on obtient
$$v_{\qq}(c)=0.$$
Pour la m\^eme raison, on a 
$$v_{\qq}(b)=0.$$
L'\'egalit\'e $2dz=b^p+c^p-a^p$ (lemme 2) et l'hypoth\`ese selon laquelle $q$ ne divise pas $W_n$, impliquent alors (comme dans la d\'emonstration du lemme 3)
$$v_{\qq}(a)\geq 1.$$
Il en r\'esulte que l'on a  
$$v_{\qq}(x+y)\geq 1.$$
On a ainsi l'\'egalit\'e
$$d(-z^p)=a^ps \quad \hbox{avec}\quad s=\sum_{k=0}^{p-1} x^{p-1-k}(-y)^k,$$
et la congruence
$$s\equiv px^{p-1}\mod \qq.$$
On en d\'eduit que  $v_{\qq}(s)=0$, donc  ${z\over a}$ est une unit\'e modulo $\qq$.  On a $dx\equiv b^p \mod \qq$, d'o\`u 
$$p\equiv {s\over x^{p-1}}\equiv \biggl({-bz\over xa}\biggr)^p \mod \qq.$$
Le corps $O_K/\qq$ est de cardinal $q$ et on a $np=q-1$, d'o\`u
$p^n\equiv 1 \mod \qq$,
puis
$$p^n\equiv 1 \mod q.$$
Puisque $n$ est pair, on  obtient
$$1=(-1)^n=(np-q)^n\equiv n^np^n\equiv n^n\mod q,$$
ce qui conduit \`a une contradiction, d'o\`u le th\'eor\`eme.
\bigskip

{\douze III. D\'emonstration du th\'eor\`eme 2}
\medskip
 Elle est analogue \`a celle du th\'eor\`eme 1 de [Ha-Pa]. Soit $\wp$ un id\'eal premier de $O_K$ au-dessus de $p$ de degr\'e r\'esiduel $1$. Notons $v_{\wp}$ la valuation sur $K$ qui lui est associ\'ee. Posons  $e=v_{\wp}(p)$ l'indice de ramification de $\wp$ sur $p$. 
Supposons qu'il existe 
$x,y,z$ dans $O_K$ tels que l'on ait
$$x^p+y^p+z^p=0\quad \hbox{et}\quad  v_{\wp}(xyz)=0.$$
Les corps $O_K/\wp$ et $\F_p$ \'etant isomorphes, il existe   $x_0,y_0,z_0\in \Z$ non divisibles par $p$ 
tels que
$$x\equiv x_0\mod \wp,\ y\equiv y_0\mod \wp, \ z\equiv z_0\mod \wp.$$
Il en r\'esulte que l'on a 
$$v_{\wp}(x^p-x_0^p)\geq \Inf( e+1,p)=e+1,$$
d'o\`u  la congruence
$$x_0^p+y_0^p+z_0^p\equiv 0 \mod \wp^{e+1}.$$
On en d\'eduit que l'on a
$$x_0^p+y_0^p+z_0^p\equiv 0 \mod p^2.$$
En particulier, on a
$$x_0+y_0+z_0\equiv 0 \mod p,$$
d'o\`u 
$$(x_0+y_0)^p+z_0^p\equiv 0 \mod p^2,$$
puis 
$$(x_0+y_0)^p\equiv x_0^p+y_0^p\mod p^2.$$
On obtient
$$1+a^p\equiv (1+a)^p\mod p^2 \quad \hbox{avec} \quad a\equiv x_0^{-1}y_0\mod p^2.\leqno(4)$$
On a $a\not\equiv 0 \mod p$ et $a\not\equiv -1 \mod p$ car $p$ ne divise pas $z_0$. Parce que (4) ne d\'epend que de la congruence de $a$ modulo $p$, on peut supposer que l'on a 
$1\leq a\leq p-2.$
 Si $a={p-1\over 2}$,  alors $a=1$ est aussi solution de  (4). Si l'on a
${p-1\over 2}<a\leq p-2$,
 alors 
 $p-1-a$ satisfait (4) et $1\leq p-1-a\leq {p-3\over 2}$. Cela contredit la condition (1), d'o\`u le r\'esultat.
\vfill\eject

{\douze IV. D\'emonstration du corollaire 3}
\medskip

1) On a $p\equiv 2 \mod 3$  car $p$ satisfait la condition (1). Dans l'anneau d'entiers de 
$\Q\bigl(\root 3\of {d}\bigr)$, il existe donc un id\'eal premier au-dessus de $p$ de degr\'e r\'esiduel 1  ([Co-1], cor. 6.4.15 et th. 6.4.16), d'o\`u la premi\`ere  assertion.
\smallskip
2) Posons $K=\Q\bigl(\root n\of {d}\bigr)$.  Le polyn\^ome  $X^n-d$ est irr\'eductible sur $\Q$, de  discriminant 
$$(-1)^{{n(n-1)\over 2}} n^nd^{n-1}.$$
La congruence  $n\equiv 1 \mod p-1$ implique  $d^n\equiv d \mod p$. 
Le fait que $p$ ne divise pas $dn$ entra\^\i ne alors  l'existence d'un id\'eal premier de $O_K$ au-dessus de $p$ non ramifi\'e  de degr\'e r\'esiduel $1$ ([Co-1], th. 4.8.13), d'o\`u le r\'esultat.
\smallskip
3) La derni\`ere assertion est  une   cons\'equence directe du th\'eor\`eme 2. 
 \bigskip
\bigskip
 
\centerline {\douze{Bibliographie}}
\bigskip
\vskip0pt\noindent
[Co-1] H. Cohen, A course in Computational Algebraic Number Theory, Springer-Verlag  GTM {\bf{138}}, 1993.
\smallskip
\vskip0pt\noindent
[Co-2] H. Cohen,  Number Theory Volume I : Tools and Diophantine Equations, Springer-Verlag  GTM {\bf{239}}, 2007.
\smallskip
\vskip0pt\noindent
[Fr-Si]  N. Freitas et S. Siksek, Modularity and the Fermat equation over totally real number fields, 
arXiv : 1307.3162v2  (2014), 32 pages.
\smallskip
\vskip0pt\noindent
[Ha-Pa] F. H. Hao, C.  J. Parry, The Fermat equation over quadratic fields, {\it J. Number Theory} {\bf 19} (1984),  115-130.
\smallskip
\vskip0pt\noindent
[Ja-Me] F. Jarvis et P. Meekin, The Fermat equation over $\Q\bigl(\sqrt{2}\bigr)$, {\it J. Number Theory} {\bf 109} (2004),  182-196.
\smallskip
\vskip0pt\noindent
[Pari]  C. Batut, D. Bernardi, K. Belabas, H. Cohen et M. Olivier, PARI-GP, version 2.3.3,  
Universit\'e de Bordeaux I, (2008).
\smallskip
\vskip0pt\noindent
[Ri] P. Ribenboim, Fermat's Last Theorem for Amateurs, Springer-Verlag, 1999.
\smallskip
\vskip0pt\noindent
[Wi]  A. Wiles, Modular elliptic curves and
Fermat's Last Theorem, {\it {Ann. of Math.}} {\bf{141}} (1995), 443-551.
\bigskip

\line {\hfill{2 avril 2014}}

\item{} Alain Kraus
\item{}Universit\'e de Paris VI, 
\item{} Institut de Math\'ematiques, 
\item{} 4 Place Jussieu, 75005 Paris,  
\item{} France
\medskip
\vskip0pt\noindent
\item{}e-mail : alain.kraus@imj-prg.fr
\smallskip

\bigskip

\bye